\documentclass[11pt]{article}
\usepackage{amsmath,amsfonts,amssymb}
\newtheorem{theorem}{Theorem}[section]
\newtheorem{lemma}[theorem]{Lemma}
\newtheorem{proposition}[theorem]{Proposition}
\newtheorem{coro}[theorem]{Corollary}
\begin{document}
\newcommand{\auteur}{
\vskip 2truecm
\centerline{Fran{\c c}ois Labourie}
\centerline{Topologie et Dynamique}
\centerline{Universit{\'e} Paris-Sud}
\centerline{F-91405 Orsay (Cedex)}}
\title{Fuchsian Affine Actions of Surface Groups}
\author{Fran{\c c}ois LABOURIE \thanks{L'auteur remercie l'Institut Universitaire de France.}}
\maketitle
\newcommand{\slr}{SL(2,{\mathbb R})}
\newcommand{\hyp}{{\mathbb H}^2}
\newcommand{\sk}{\section}
\newcommand{\ssk}{\subsection}
\newcommand{\sssk}{\subsubsection}
\newcommand{\nbb}{\overline{\nabla}}
\newcommand{\inn}{\}_{n\in{\mathbb N}}}
\newcommand{\bull}{\textsc q.e.d.}
\newcommand{\rr}{{\mathbb R}}
\newcommand{\cc}{{\mathbb C}}
\newcommand{\hh}{{\mathbb H}}
\newcommand{\pp}{{\mathbb P}}
\newcommand{\qq}{{\mathbb Q}}
\newcommand{\zz}{{\mathbb Z}}
\newcommand{\nn}{{\mathbb N}}
\newcommand{\vv}{{\mathbb V}}
\renewcommand{\a}{\alpha}
\renewcommand{\b}{\beta}
\newcommand{\g}{\gamma}
\renewcommand{\k}{\kappa}
\renewcommand{\l}{\lambda}
\newcommand{\Si}{\Sigma}
\newcommand{\s}{\sigma}
\newcommand{\e}{\varepsilon}
\newcommand{\D}{\Delta}
\renewcommand{\d}{\delta}
\renewcommand{\L}{\Lambda}
\newcommand{\G}{\Gamma}
\renewcommand{\r}{\rho}
\newcommand{\nb}{{\nabla}}
\newcommand{\begitem}{\begin{itemize}}
\newcommand{\enitem}{\end{itemize}}
\newcommand{\iti}{\item[(i)]}
\newcommand{\itii}{\item[(ii)]}
\newcommand{\itiii}{\item[(iii)]}
\newcommand{\itiv}{\item[(iv)]}
\newcommand{\itv}{\item[(v)]}
\newcommand{\itvi}{\item[(vi)]}
\newcommand{\itvii}{\item[(vii)]}
\newcommand{\ita}{\item[(a)]}
\newcommand{\itb}{\item[(b)]}
\newcommand{\itc}{\item[(c)]}
\newcommand{\itd}{\item[(d)]}
\newcommand{\ite}{\item[(e)]}
\newcommand{\itf}{\item[(f)]}
\newcommand{\itg}{\item[(g)]}
\newcommand{\ith}{\item[(h)]}
\newcommand{\tit}{\item[-]}
\newcommand{\preu}{{\sl Proof:~}}
\newcommand{\theo}[2]{\begin{theorem}
\label{#1}
#2
\end{theorem}}

\newcommand{\lem}[2]{
\begin{lemma}
\label{#1}
#2
\end{lemma}}

\newcommand{\cor}[2]{\begin{coro}
\label{#1}
#2
\end{coro}}

\newcommand{\pro}[2]{
\begin{proposition}
\label{#1}
#2
\end{proposition}}

\newcommand{\nin}{n\in{\mathbb N}}
\newcommand{\lla}{\lfloor}
\newcommand{\rra}{\rfloor}
\sk{Introduction}

Let $\Gamma$ be the fundamental group of a compact surface. Let $\l_{n}$ be the standard $n$-dimensional representation of $\slr$ in $SL(n,\rr)$. We shall say a representation $\rho$ of $\G$ in $SL(n,\rr)$ is {\it Fuchsian} ( or {\it $n$-Fuchsian }) if $\rho=\lambda_{n}\circ \iota$, where $\iota$ is a discrete faithful representation of $\G$ in $SL(2,\rr)$. We shall also say by extension the image of $\rho$ is {\it Fuchsian}, and that an affine action of a surface group is  {\it Fuchsian}, if its linear part is Fuchsian.

Our main result is the following theorem

\theo{maintheorem}{A finite dimensional affine  Fuchsian action of the fundamental group of a compact surface is not proper.}

In even dimensions, this is a trivial remark. For dimension $4p+1$, this theorem follows from previous results of G.~Margulis \cite{Mar1} \cite{Mar2} concerning actions of free groups by using  the {\it Margulis invariant} and  lemma \ref{lemmarg}, also due to Margulis. Therefore, our proof shall concentrate on dimensions $4p+3$.

This case bears special features: one should notice that  G.~Margulis has 
exhibited proper actions of free group (with two generators) on ${\rr}^{3}$ \cite{Mar1} \cite{Mar2}. Therefore, surface groups behave differently than free groups  in these dimensions. 

When 
$\dim(E)=3$, our result is a theorem of G.~Mess \cite{Mess}, for which G.~Margulis and W.~Goldman \cite{MarGold} have obtained a different proof using Margulis invariant and Teichm\"uller theory. Our proof is based on similar ideas, but uses instead of Teichm\"uller theory a result on Anosov flows and a holomorphic interpretation of Margulis invariant, hence generalizing to higher dimensions.

It is a pleasure to thank M.~Babillot, W.~Goldman, G.~Margulis for helpful conversations. 

\sk{Representations of $SL(2,{\mathbb R})$, surfaces and connections}

In this section, we describe the irreducible representation of $SL(2,{\mathbb R})$ of dimension $2n+1$ as the holonomy  of a flat connection.

It is well known that in dimension 3, the irreducible representation of $\slr$
is associated with the Minkowski model of the hyperbolic plane $\hyp$. More precisely, there exists a flat connection on 
$E=\rr\oplus T\hyp$, such that the action of $\slr$ lifts to a connection preserving action on this bundle. Hence, we obtain  a 3-dimensional representation of $\slr$. Furthermore, the Minkowski model is obtained using the section $(1,0)$ of $E$.

We will now be more precise and explain this construction in more details in higher dimensions.

\ssk{A flat connection}{\label{connection}}

Let $\hyp$ be the hyperbolic plane with its complex structure.
Let $L_{k}$ be the complex line bundle over $\hyp$ defined by 
$$
L_{k}=(T\hyp)^{\otimes_{\mathbb C}^{k}}.
$$ 
Let 
$$
E={\mathbb R}\oplus L_{1}\oplus\ldots\oplus L_{n}.
$$
Notice now that $\slr$ acts on all $L_{k}$, hence on $E$,  by bundle automorphisms.

If $Y$ is a section of $E$, $Y_{0}$ will denote its component on the factor $\rr$, and $Y_{k}$ its component on $L_{k}$. The space of sections of the  bundle $V$ will be denoted $\Gamma(V)$.
The metric on $L_{i}$, induced from the Riemanian metric on $\hyp$ will be denoted $\langle,\rangle$. By definition, if $Y\in L_{k}$, $X\in L_{1}$, then $i_{X}Y$ is the element of $L_{k-1}$ such that
$$
\forall Z\in L_{k},~ ~ \langle i_{X}Y,Z\rangle=\langle Y,X\otimes Z\rangle.
$$
Let $\nbb$ be the Levi-Civita connection on $L_{1}$, and, by extension, the induced connection on $L_{k}$. We introduce the following connection $\nb$ on $E$, defined if $X\in T\hyp$, $Y\in \Gamma(E)$ by 
\begin{displaymath}
\left\{
\begin{array}{rcl}
(\nb_{X}Y)_{0}&=&L_{X}Y_{0}+\frac{1}{2}(n+1)\langle X,Y_{1}\rangle\\
\forall k>0,\,\, (\nb_{X}Y)_{k}&=&(n-k+1)X\otimes Y_{k-1} + \nbb_{X}Y_{k}\\
& & + \frac{1}{4}(n+k+1)i_{X}Y_{k+1}.
\end{array}
\right.
\end{displaymath}
Consider the family or real numbers, defined for  $k\in\{0,n-1\}$,
$$
a_{0}=1,~a_{k+1}=\frac{1}{2^{2k+1}}\prod_{j=0}^{j=k}(\frac{n+j+1}{n-j}).
$$
Define a metric of signature $(n,n+1)$ on $E$ by
$$
\lla Y,Z\rra=\sum_{k=0}^{k=n}(-1)^{k+1}a_{k}\langle Y_{k},Z_{k}\rangle.
$$
The main result of this section is the following statement
\pro{flatconnexion}{The connection $\nb$ is flat, and preserves the metric $\lla,\rra$. Furthermore, the $\slr$ action on $E$ preserves the metric $\lla,\rra$ and the connection $\nb$. The resulting $(2n+1)$-representation of $\slr$ is irreducible.}
\preu Let $X$, $Z$ two commuting vector fields on $\hyp$. Let $\omega$ the K\"ahler form of $\hyp$ defined by $\omega(Z,X)=\langle JZ,X\rangle$. Let's first introduce the following notation. If $f$ is a function of $Z$ and $X$ then
\newcommand{\www}{\underline}
$$
\www{f(Z,X)}=f(Z,X)-f(X,Z).
$$
With these notations at hands, we have
\begin{eqnarray*}
\www{Z\otimes\langle X,Y\rangle}&=&-\omega(Z,X)JY,\\
\www{Z\otimes i_{X}Y}&=&-2\omega(Z,X)JY,\\
\www{i_{Z}(X\otimes Y)}&=&2\omega(Z,X)JY.
\end{eqnarray*}
Let $Y$ be a section of $E$. Let $\bar{R}$ be the curvature tensor of $\nbb$ and recall
that
$$
\bar{R}(Z,X)Y_{k}=k\omega(Z,X)JY.
$$
We begin our computations. Let $R$ be the curvature tensor of $\nb$.
We first have
\begin{eqnarray*}
\lefteqn{(\nb_{Z}\nb_{X}Y)_{0}=}\\
& & L_{Z}L_{X}Y_{0}+\frac{1}{2}(n+1)L_{Z}\langle X,Y_{1}\rangle\\
& & +\frac{1}{2}(n+1)(n)\langle Z,X\rangle \otimes Y_{0}+ \frac{1}{2}(n+1)\langle Z,\nbb_{X}Y_{1}\rangle \\
& & +\frac{1}{8}(n+1)(n+2)\langle Z,i_{X}Y_{2}\rangle.
\end{eqnarray*}
Hence
$$
(R(Z,X)Y)_{0}=0.
$$
Next
\begin{eqnarray*}
\lefteqn{(\nb_{Z}\nb_{X}Y)_{1}=}\\
& & nZ\otimes L_{X}Y_{0} +\frac{1}{2}n(n+1)Z\otimes\langle X,Y_{1}\rangle
+n\nbb_{Z}(X\otimes Y_{0})\\
& & +\nbb_{Z}\nbb_{X}Y_{1}+\frac{1}{4}(n+2)\nbb_{Z}(i_{X}Y_{2})
+\frac{1}{4}(n+2)(n-1)i_{Z}(X\otimes Y_{1})\\
& &+\frac{1}{4}(n+2)i_{Z}\nbb_{X}Y_{2}+\frac{1}{16}(n+2)(n+3)i_{Z}i_{X}Y_{3}.
\end{eqnarray*}
We get
\begin{eqnarray*}
\lefteqn{(R(Z,X)Y)_{1}=}\\
& & +\frac{1}{2}n(n+1)\www{Z\otimes\langle X,Y_{1}\rangle} +
\bar R(Z,X)Y_{1}\\
& & +\frac{1}{4}(n+2)(n-1)\www{i_{Z}(X\otimes Y_{1})}
+\frac{1}{16}(n+2)(n+3)\www{i_{Z}i_{X}Y_{3}}.
\end{eqnarray*}
Hence
\begin{eqnarray*}
\lefteqn{(R(Z,X)Y)_{1}}\\
& & =\omega(Z,X)JY\big(-\frac{1}{2}n(n+1)+1 +\frac{2}{4}(n+2)(n-1)\big)\\
& & =0.
\end{eqnarray*}
It remains to consider the case $k>0$
\begin{eqnarray*}
\lefteqn{(\nb_{Z}\nb_{X}Y)_{k}=}\\
& & (n-k+1)Z\otimes \nbb_{X}Y_{k-1} +\frac{1}{4}(n-k+1)(n+k)Z\otimes i_{X}Y_{k}\\
& & +(n-k+1)(n-k+2)Z\otimes X\otimes Y_{k-2}
+(n-k+1)\nbb_{Z}(X\otimes Y_{k-1})\\
& & +\nbb_{Z}\nbb_{X}Y_{k}+\frac{1}{4}(n+k+1)\nbb_{Z}(i_{X}Y_{k+1})\\
& &+\frac{1}{4}(n+k+1)(n-k)i_{Z}(X\otimes Y_{k})+\frac{1}{4}(n+k+1)i_{Z}\nbb_{X}Y_{k+1}\\
& & +\frac{1}{16}(n+k+1)(n+k+2)i_{Z}i_{X}Y_{k+2}.
\end{eqnarray*}
We get
\begin{eqnarray*}
\lefteqn{(R(Z,X)Y)_{k}=}\\
& & +\frac{1}{4}(n-k+1)(n+k)\www{Z\otimes i_{X}Y_{k}} +
\bar R(Z,X)Y_{k}\\
& & +\frac{1}{4}(n+k+1)(n-k)\www{i_{Z}(X\otimes Y_{k})}\\
& & +\frac{1}{16}(n+k+1)(n+k+2)\www{i_{Z}i_{X}Y_{k+2}}\\
& & +(n-k+1)(n-k+2)\www{(Z\otimes X\otimes Y_{k-2})}.
\end{eqnarray*}
Hence
\begin{eqnarray*}
\lefteqn{(R(Z,X)Y)_{k}}\\
& & =\omega(Z,X)JY\big(-\frac{2}{4}(n-k+1)(n+k)+k +\frac{2}{4}(n+k+1)(n-k)\big)\\
& & =0.
\end{eqnarray*}
We have just proved the connection $\nbb$ is flat.
Now, we show $\nb$ preserves $\lla,\rra$. Let $Y$ a section of $E$. Then
\begin{eqnarray*}
\lefteqn{\lla \nb_{X}Y,Y\rra}\\
& = & \sum_{k=0}^{k=n}(-1)^{k+1}a_{k}\langle (\nb_{X})Y_{k},Y_{k}\rangle\\
& = & -\langle L_{X}Y_{0}-\frac{1}{2}(n+1)\langle X,Y_{1}\rangle,Y_{0}\rangle+\sum_{k=1}^{k=n}(-1)^{k+1}a_{k}\langle \nbb_{X}Y_{k},Y_{k}\rangle\\
& & +\sum_{k=1}^{k=n}(-1)^{k+1}a_{k}\langle (n-k+1)X\otimes Y_{k-1} +  
\frac{(n+k+1)}{4}i_{X}Y_{k+1},Y_{k}\rangle\\
& = & -\langle L_{X}Y_{0},Y_{0}\rangle 
+\sum_{k=1}^{k=n}(-1)^{k+1}a_{k}\langle \nbb_{X}Y_{k},Y_{k}\rangle\\
& &-\frac{1}{2}(n+1)\langle \langle X,Y_{1}\rangle,Y_{0}\rangle \\
& & +\sum_{k=1}^{k=n}(-1)^{k+1}a_{k}\frac{(n+k+1)}{4}
\langle i_{X}Y_{k+1},Y_{k}\rangle\\
& &+\sum_{k=1}^{k=n}(-1)^{k+1}a_{k}(n-k+1)\langle X\otimes Y_{k-1},Y_{k}\rangle.
\end{eqnarray*}
We make a change of variables in the very last sum, and get
\begin{eqnarray*}
\lefteqn{\lla \nb_{X}Y,Y\rra}\\
&= &-\langle L_{X}Y_{0},Y_{0}\rangle 
+\sum_{k=1}^{k=n}(-1)^{k+1}a_{k}
\langle \nbb_{X}Y_{k},Y_{k}\rangle\\
& & -\frac{1}{2}(n+1)\langle \langle X,Y_{1}\rangle,Y_{0}\rangle \\
 & & +\sum_{k=1}^{k=n}(-1)^{k+1}a_{k}\frac{(n+k+1)}{4}
\langle i_{X}Y_{k+1},Y_{k}\rangle\\
& & +\sum_{k=0}^{k=n-1}(-1)^{k}a_{k+1}(n-k)\langle X\otimes Y_{k},Y_{k+1}\rangle
\end{eqnarray*}
Hence
\begin{eqnarray*}
\lefteqn{\lla \nb_{X}Y,Y\rra}\\
& =&  L_{X}\lla Y,Y\rra\\
& & -\frac{1}{2}(n+1)\langle \langle X,Y_{1}\rangle,Y_{0}\rangle +n a_{1}\langle X\otimes Y_{0},Y_{1}\rangle\\
& & \sum_{k=1}^{k=n}(-1)^{k}\big(a_{k+1}(n-k)-a_{k}\frac{(n+k+1)}{4}\big)\langle X\otimes Y_{k},Y_{k+1}\rangle.
\end{eqnarray*}
To conclude, we just have to remark that
\begin{eqnarray*}
a_{1}&=&\frac{n+1}{2n}\\
\frac{a_{k+1}}{a_{k}}&=&\frac{n+k+1}{4(n-k)}.
\end{eqnarray*}
We finally have to check that the corresponding representation of the group $\slr$ is irreducible. For that let $S^{1}\subset \slr$, a subgroup isomorphic to the circle fixing a point
$x_{0}$. The corresponding action on $L_{k}(x_{0})$ is given by
$$
e^{i\theta}(u)=e^{k i\theta}u.
$$
This shows the representation is the $2n+1$ dimensional one\bull

\sk{Cohomology and $(n+1)$-holomorphic differentials }
Let $S=\hyp/\Gamma$ be a compact surface.  Let  $\r$ be a $(2n+1)$-Fuchsian representation of $\Gamma$. In this section,  we shall describe  the vector space $H^{1}_{\r}(\Gamma,{\mathbb R}^{2n+1})$ in terms of $(n+1)$-holomorphic differentials on $S$.

We use the notations of the previous sections. Let $E_{S}=E/\Gamma$ be the vector bundle over $S=\hyp/\Gamma$ coming from $E$.

Let ${\mathcal H}^{q}$ the vector space of $q$-holomorphic differentials on $S$. Let $\Lambda^{p}(E_{\Gamma})$ the vector space of $p$-forms on $S$ with value in $E_{S}$. The flat connection $\nb$ gives rise to a complex
$$
0\longrightarrow\Lambda^{0}(E)\stackrel{d^{\nb}}{\longrightarrow}\Lambda^{1}(E)\stackrel{d^{\nb}}{\longrightarrow}\Lambda^{2}(E)\longrightarrow 0.
$$
The cohomology of this complex is $H^{*}_{\rho}(\Gamma,{\mathbb R}^{n+1})$.
From the metric on $\hyp$, we deduce an isomorphism $\omega\mapsto\check\omega$ of  $L_{k}^{*}$ with $L_{k}$. We define  now a map $\Phi$ by
\begin{displaymath}
\Phi : \left\{
\begin{array}{rcl}
{\mathcal H}^{2n+1}&\rightarrow&\Lambda^{1}(E)\\
\omega&\mapsto&(X\mapsto i_{X}\check{\omega}\in L_{n}\in E)
\end{array}
\right.
\end{displaymath}
We first prove:
\pro{injectivity}{For every $(n+1)$-holomorphic differential
$\omega$
$$
d^{\nabla}(\Phi(\omega))=0.
$$
Furthermore, if $\Phi(\omega)=d^{\nabla}u$, then $\omega=0$.}
\preu
By definition,
$$
d^{\nabla}\Phi(\omega)(X,Y)=\nabla_{X}i_{Y}\check{\omega}-\nabla_{Y}i_{X}\check{\omega}-i_{[X,Y]}\check{\omega}.
$$
Hence, if $n>1$
$$
d^{\nabla}\Phi(\omega)(X,Y)
=\frac{2n}{4}(i_{X}i_{Y}\check{\omega}-i_{Y}i_{X}\check{\omega})+(
\nbb_{X}i_{Y}\check{\omega}-\nbb_{Y}i_{X}\check{\omega}-i_{[X,Y]}\check{\omega}).
$$
Notice that $i_{X}i_{Y}\check{\omega}$ is symmetric in $X$ and $Y$. Finally, the holomorphicity condition on $\omega$ implies 
$$
\nbb_{X}i_{Y}\check{\omega}-\nbb_{Y}i_{X}\check{\omega}-i_{[X,Y]}\check{\omega}=0.
$$
A similar proof (but with different constants) yields the result for $n=1$.

Next, assume $\Phi(\omega)=d^{\nb}u$. The (non Riemaniann) metric on $E$ and the Riemannian metric on $\hyp$ induce a metric on $\Lambda^{*}(E)$, which we denote $\lla,\rra_{\Lambda}$. One should notice here that even though this metric is neither positive nor negative, since $\Phi(\omega)$ is a section of a bundle on which the metric is either positive or negative, we have
$$
\lla \Phi(\omega), \Phi(\omega)\rra_{\Lambda}=0
\Rightarrow\Phi(\omega)=0
\Rightarrow\omega=0
$$  
Let $(d^{\nb})^{*}$ be the adjoint of $d^{\nb}$. One has, if $(X_{1},X_{2})$ is a basis of $T\hyp$,
$$
(d^{\nabla})^{*}(\phi(\omega))=-\sum_{k=1}^{k=2}\nabla_{X_{k}}(i_{X_{k}}\check{\omega}).
$$
A short calculation shows
$$
(d^{\nabla})^{*}(\phi(\omega))=-\sum_{k=1}^{k=2}\nbb_{X_{k}}(i_{X_{k}}\check{\omega}),
$$
and this last term is $0$ by holomorphicity.
We have just proved that 
$$
(d^{\nabla})^{*}\Phi(\omega)=0.
$$ 
Hence, $\Phi(\omega)=d^{\nb}u$ implies 
$$
\lla\Phi(\omega), \Phi(\omega)\rra_{\Lambda}=
\lla (d^{\nabla})^{*}\Phi(\omega),u\rra_{\Lambda}=0.
$$
This ends the proof \bull

It follows from the previous proposition that $\Phi$ gives rise to a map (also denoted $\Phi$) from ${\cal H}^{n+1}$ to the space $H^{1}_{\r}(\Gamma, \rr^{2n+1})$. We have:

\cor{isocoho}{The map $\Phi$ is an isomorphism from ${\cal H}^{n+1}$ to the space  $H^{1}_{\r}(\Gamma, \rr^{2n+1})$.}

\preu Indeed, we have just proved that $\Phi$ is injective. Furthermore, if  $\chi(S)$ is the Euler characteristic of $S$, we have 
$$
\dim(H^{1}_{\r}(\Gamma, \rr^{2n+1}))\leq (2n+1)\chi(S).
$$
But, by Riemann-Roch, 
$$
\dim({\cal H}^{n+1})=(2n+1)\chi(S).
$$
Hence, the corollary follows \bull

\sk{A de Rham interpretation of Margulis invariant}

The irreducible representation  of $\slr$ of dimension  $2n+1$ preserves a metric $\lla,\rra$ of signature $(n, n+1)$.
\ssk{Loxodromic elements}
We define a {\it loxodromic}  element in $SO(n,n+1)$ to be ${\mathbb R}$-split and in the interior of a Weyl chamber. This just means all eigenvalues are real and  have multiplicity 1. Recall that 1 allways belong to the spectrum of a loxodromic element. Notice that all the elements, except the identity, of a $(2n+1)$-Fuchsian surface group are loxodromic.

\ssk{The invariant vector of a loxodromic element}
Chose now, once and for all,  an orientation on ${\mathbb R}^{2n+1}$. The light cone - without the origin - has two components. Let's also choose one of these components.

Let $\gamma$ be  a loxodromic element.
It follows from the previous choices  that we have a well defined eigenvector, the {\it invariant vector}, denoted  $v_{\gamma}$,
associated to the eigenvalue 1.

Indeed, all the other eigenvectors are lightlike. We order  all the eigenvaluesdifferent than 1, in such a way that $\l_{i}<\l_{i+1}$.
Thanks to our choices, we may pick one eigenvector $e_{i}$ 
in the preferred component of the light cone for all the eigenvalues $\l_{i}$ different than 1. We now choose $v_{\gamma}$ of norm 1, such that $(v_{\gamma},e_{1},\ldots,e_{2n})$ is positively oriented.

\ssk{Margulis invariant} 
Let now $Iso(n,n+1)={\mathbb R}^{2n+1}\rtimes SO(n,n+1)$ be 
the group of orientation preserving isometries of  ${\mathbb R}^{2n+1}$ as an affine space. For $\gamma$ in $Iso(n,n+1)$,
$\hat\gamma$ denotes  its linear part. 
We shall say an element of $Iso(n,n+1)$ is {\it loxodromic} if its linear part is a loxodromic element of $SO(n,n+1)$.

The {\it Margulis invariant} ( \cite{Mar1} \cite{Mar2} ) of a loxodromic element $\gamma$ of $Iso(n,n+1)$ is
$$
\mu(\gamma)=\lla\gamma (x)-x,v_{\hat\gamma}\rra,
$$
where $x$ is an element of ${\mathbb R}^{2n+1}$. A quick check shows 
$\mu(\gamma)$ does not depend on $x$. 

\ssk{Margulis invariant and properness of an affine action}

Let $\gamma_{1}$ and $\gamma_{2}$ be two loxodromic elements.
Let $E^{+}_{i}$ (resp. $E^{-}_{i}$)  be the space generated by the eigenvectors of $\hat\gamma_{i}$ corresponding to  the eigenvalues of absolute value greater (resp. less) than 1.

We say  $\gamma_{1}$ and $\gamma_{2}$ are in {\it general position} if the two decompositions 
$$
{\mathbb R}.v_{\hat\gamma_{i}}\oplus E^{+}_{i}\oplus E^{-}_{i},
$$
are in general position.

Notice that for a $(2n+1)$-Fuchsian group, two (non comensurable) elements are loxodromic and in general position.

In \cite{Mar1} \cite{Mar2} (see also \cite{Drumm}) G. Margulis has proved the following magic lemma 

\lem{lemmarg}{If two loxodromic elements $\gamma_{1}$, $\gamma_{2}$, in general position, are such that $\mu(\gamma_{1})\mu(\gamma_{2})\leq 0$,  then the group generated by
$\gamma_{1}$ and $\gamma_{2}$ does not act properly on ${\rr}^{2n+1}$} 

\ssk{An interpretation of Margulis invariant}

Let $\r$ a representation of $\Gamma$ in $Iso(n,n+1)$, whose linear part, $\hat\rho$, is Fuchsian.
Let $E_{S}={\mathbb R}\oplus L_{1}\oplus\ldots\oplus L_{n}$ the flat bundle over $S$ described in \ref{flatconnexion}  whose holonomy is  $\hat\rho$. 

We describe now  $\rho$ as  an element of $H^{1}_{\hat\r}(\Gamma,\rr^{2n+1})$. 

Let  $\alpha\in H^{1}_{\r}(\Gamma,\rr^{2n+1})$, interpreted as an element of $\Lambda^{1}(E_{S})$. Let $\nb^{\alpha}$ be the flat connection on 
$F={\mathbb R}\oplus E_{S}$ defined by
$$
\nb^{\alpha}_{X}(\lambda, V)=(L_{X}\lambda, \lambda.\alpha(X)+\nb_{X}V).
$$
We claim there exists  $\alpha\in H^{1}_{\hat\r}(\Gamma,\rr^{2n+1})$ such that the holonomy of $\nb^{\alpha}$ is $\rho$. Of course, here,  $\rr^{p}\rtimes SL(p,\rr)$ is identified with a subgroup of $GL(p+1,\rr)$.

Let now  $c$ be a closed curve on $S$, represented in homotopy by the conjugacy class of some element $\gamma$. Since $v_{{\hat\r}(\gamma)}$ is invariant under ${\hat\r}(\gamma)$, it gives rise to a parallel section $v_{c}$ of $E\vert_{c}$.

We first prove the following statement:

\pro{inter1}{Let $c$, $\r$, $\gamma$, $\a$ be as above. Then
$$
\mu(\r(\gamma))=\int_{c}\lla\a,v_{c}\rra.
$$
}

\preu We shall use the previous notations. We parametrise $c$ by the circle of length 1. Let $\pi$ be the covering $\hyp\rightarrow S$. Consider a lift $\tilde c$ of $c$ on the universal cover of $S$. The bundle $\pi^{*}F$ becomes trivial. The canonical section $\sigma$ corresponding to the $\rr$ factor in $F$, gives rise to a map 
$$
i:\hyp\rightarrow {\rr}^{2n+1},
$$
taking value in the affine hyperplane 
$$
P=\{(1,u)\in{\rr}^{2n+1}\}.
$$
Let ${\bar c}=i\circ{\tilde c}$, and let's identify $\r(\g)$ with $\g$. 
By definition now:
\begin{eqnarray*}
\mu(\gamma)&=&\lla{\r(\g)(\bar c}(0))-{\bar c}(0),v_{\hat\gamma}\rra\\
& =  &\lla{\bar c}(1)-{\bar c}(0),v_{\hat\gamma}\rra\\
& =  &\int_{0}^{1}\lla{\dot{\bar c}}(s),v_{\hat\gamma}\rra ds
\end{eqnarray*}
Now, we interpret the last term  on $F$ and we obtain 
\begin{eqnarray*}
\mu(\gamma)& = &\int_{0}^{1}\lla\nb^{\a}_{{\dot c}(s)}\sigma,v_{c}(s)\rra ds\\
& =  & \int_{0}^{1}\lla\a({\dot c}(s)),v_{c}(s)\rra ds\\
& =  & \int_{c}\lla\a,v_{c}\rra.
\end{eqnarray*}
This ends the proof \bull

\ssk{The invariant vector as a section}

In this paragraph, we assume $n=2p+1$, so that our representation is of dimension $4p+3$.  

We use the notations of the previous paragraphs. In particular, let $\g\in\Gamma$. Let $v=v_{{\hat\r}({\gamma})}$. Let $c$ be the closed geodesic (for the hyperbolic metric) corresponding to the element $\g$. 

Recall that $v_{\gamma}$ gives rise to a section $v_{c}$ along the closed geodesic, which is parallel. 

In this paragraph, we wish to describe $v_{c}$ explicitely. Let $J$ the complex structure of $S$.  Let's  introduce the following section (along c) defined by

\begin{eqnarray*}
(w_{c})_{2k} & = & 0 \\
(w_{c})_{2k+1} & = & 
J(-4)^{k}\prod_{l=1}^{l=k}(\frac{p-l}{p+l+1}) \,\,\underset{2k+1}{\underbrace{{\dot c}\otimes\ldots\otimes{\dot c} }}.
\end{eqnarray*} 

\pro{sectiongeodesic}{The section $w_{c}$ of $E_{S}$  is parallel along $c$. Furthermore, there exists $\varepsilon\in\{-1,1\}$ independant of $c$ such that
$$
v_{c}=\varepsilon\frac{w_{c}}{\sqrt{\lla w_{c},w_{c}\rra}}.
$$}
\preu A straightforawrd computation shows that $w_{c}$ (hence $v_{c}$) is parallel. Furthermore $w_{c}$ is a space like vector, and by construction $v_{c}$ has norm 1. 

It remains to prove that $v_{c}$ has the correct orientation. 
For that consider any geodesic arc $u$ in $\hyp$ paramatrised by $[0,L]$. 
We have  a basis  of $E\vert_{u(t)}$ given by
$$
B(t)=(1,{\dot u},J{\dot u},\ldots,
\underset{n}{\underbrace{{\dot u}\otimes\ldots\otimes{\dot u}}}, J\underset{n}{\underbrace{{\dot u}\otimes\ldots\otimes{\dot u}}}).
$$
We may now consider the isometry $\gamma(u)$ sending $B(0)$ to $B(L)$. This is a loxodromic isometry. Next, consider the following section of $E$ along $u$ given by
\begin{eqnarray*}
(w_{u})_{2k} & = & 0 \\
(w_{u})_{2k+1} & = & 
J(-4)^{k}\prod_{l=1}^{l=k}(\frac{p-l}{p+l+1}) \,\,\underset{2k+1}{\underbrace{{\dot u}\otimes\ldots\otimes{\dot u} }}.
\end{eqnarray*}  

This section is parallel along $u$ and therefore gives rise to a vector proportional to the invariant vector of $\gamma(u)$.

Next, by continuity, this  proportional is constant. Applying this remark to a lift in the universal cover of our closed geodesic,  this ends the proof \bull

\sk{Main theorem in dimension $4p+3$}
Again, let's $\r$ be a representation of a compact surface group $\G$ in the group of affine transformation of an affine space of dimension $4p+3$, whose linear part $\hat\r$ is Fuchsian. 
We assume that $\r(\G)$ acts properly on  $\rr^{4p+3}$.
The representation $\r$ is described  from $\hat\r$ as an element $\a$ of
$H^{1}_{\hat\r}(\G,\rr^{4p+3})$.

According to proposition \ref{isocoho}, this element $\a$ is described by a $(2p+2)$-holomorphic differential $\omega$.

Let $\g\in \G$, and $c$ the corresponding closed geodesic. 
From proposition \ref{inter1}, we get
$$
\mu(\r(\g))=\int_{c}\lla \a, v_{c}\rra.
$$

From proposition \ref{sectiongeodesic}, we deduce there exist a constant $K_{1}$ just depending on $p$ such that
$$
\mu(\r(\g))= K_{1} \int_{c}\lla i_{\dot c}\omega, J \underset{2p+1}{\underbrace{{\dot c}\otimes\ldots\otimes{\dot c} }} \rra dt.
$$

From the constructions explained in the paragraph \ref{connection}, we finally obtain there exist a contant $K_{2}$ just depending on $p$  such that
\begin{eqnarray*}
\mu(\r(\g))&= &K_{2}\int_{c}\langle i_{\dot c}{\check\omega},J \underset{2p+1}{\underbrace{{\dot c}\otimes\ldots\otimes{\dot c} }} \rangle dt\\
&=&-K_{2}\int_{c}\Im (\omega(\underset{2p+2}{\underbrace{{\dot c}\otimes\ldots\otimes{\dot c} }}))dt.
\end{eqnarray*}

Let $US$ be the unit tangent bundle of $S$. Let $f$ be the function defined on $US$
by
$$
f(u)=\Im (\omega(\underset{2p+2}{\underbrace{{u}\otimes\ldots\otimes{u} }})).
$$

From lemma \ref{lemmarg} and the previous computation, we obtain that  the integral of $f$ along closed orbits of the geodesic flow has a constant sign. On the other hand, let $\lambda$ be the Lebesgue measure,  we have

$$
\int_{US}f d\lambda =0.
$$
Indeed, let $\beta$ be a complex number such that $\beta^{2p+2}=-1$. This number $\beta$ is to  be considered as a diffeomorphism of $US$, which preserves the orientation and the Lebesgue measure. Lastly $ f\circ \beta=-f$ and this proves the last formula.

The conclusion of the proof follows at once from the following lemma.

\lem{classic}{Let $M$ be a compact manifold equipped with an  Anosov flow preserving a measure $\nu$ which charges open sets. 
Let $f$ be a H\"older function defined on $M$ such that its integral on every closed orbit is positive, then the integral of $f$ with respect to $\nu$ is positive.}

\preu I could not find a proper reference in the litterature of this specific lemma, although lots of versions exist for discrete time transformations. Let sketch a proof by overkilling using hints from a conversation with G. Margulis. 
Let $\phi_{t}$ be the flow.
It follows from a theorem of M. Ratner \cite{Ratner}
that if $f$ is a H\"older function whose integral is 0, either it is a cocycle (and in particular, its integral over every closed orbit is 0), or it satisfies the central limit theorem. In particular, this implies that there exists at leaat one point $x$ with a dense orbit, such that

$$
\lim_{t\rightarrow\infty}\frac{1}{\sqrt t}\int_{0}^{t}f\circ\phi_{s}(x) ds<A<0.
$$

Let now $\{t_{n}\}_{n\in {\mathbb N}}$ be a sequence of real numbers converging to infinity such that  $\{\phi_{t_{n}}(x)\}_{n\in {\mathbb N}}$ converges to $x$. By the closing lemma and classical estimates, there exists $n_{0}$ such that for $n>n_{0}$ we can find a closed geodesic $c_{n}$ such that
$$
\vert\int_{0}^{t_{n}}f\circ\phi_{s}(x)ds - \int_{c_{n}}f ds\vert < A/2.
$$
This implies the lemma.\bull

\sk{Other dimensions}
The other dimensions are either trivial (even case) or follows from the immediate use of Margulis invariant ($4p+1$ case) as we shall explain now.

Let $\lambda$ be the representation of $\slr$ of even dimension. Let $h$ be a loxodromic element of $SL(2,\rr)$. Then $1$ will not belong to the spectrum of $\lambda(h)$. It follows, that if  $\rho$ is a representation of $\Gamma$ in even dimension  whose linear part is Fuchsian then for all $\g$ in $\G$ different than the identity then $\rho(\g)$ does not act properly.

Last, in dimensions $4p+1$, the Margulis invariant is such that 
$\mu(\g^{-1})=-\mu(\g)$. It follows at once from lemma \ref{lemmarg}, that if  $\rho$ is a representation of $\Gamma$ in dimension $4p+1$  whose linear part is Fuchsian, if  $\g_{1}$ and $\g_{2}$ are non commensurable elements of $\G$ then $\r(\g_{1})$ and $\r(\g_{2})$ generate a group that does not act properly on the affine space. Of course, the point in our previous discussion in that in dimension $4p+3$ then $\mu(\g^{-1})=\mu(\g)$, hence such an argument do not work and actually, free groups (even Fuchsian ones) can act properly, see \cite{Mar1}, \cite{Mar2} and \cite{Drumm}.

\auteur
\end{document}